\documentclass[a4paper,10pt]{article}
\usepackage{paper-en}
\usepackage{hyperref}

\def\thetitle{All-set-homogeneous spaces}
\def\theauthors{Nina Lebedeva and Anton Petrunin}

\hypersetup{colorlinks=true,
citecolor=black,
linkcolor=black,
anchorcolor=black,
filecolor=black,
menucolor=black,
urlcolor=black,
pdftitle={\thetitle},
pdfauthor={\theauthors}
}

\begin{document}

\title{\thetitle}
\author{\theauthors}
\date{}
\maketitle

\begin{abstract}
A metric space is said to be all-set-homogeneous if any isometry between its subsets can be extended to an isometry of the whole space.
We give a classification of a certain subclass of all-set-homogeneous length spaces.
\end{abstract}

\section{Main result}

The distance between two points $x$ and $y$ in a metric space $M$ will be denoted by $|x-y|_M$.
Recall that $M$ is called \emph{length} (or \emph{geodesic}) space if any two points $x,y\in M$ can be connected by a path $\gamma$ 
such that $|x-y|_M$ is arbitrarily close to the length of $\gamma$ (or $|x-y|_M=\length\gamma$ respectively).
Evidently, any geodesic space is a length space, but not the other way around.

A metric space $M$ is said to be \emph{all-set-homogeneous} 
if any isometry $A\to A'$ between its subsets can be extended to an isometry of the whole space $M\to M$.

Examples of geodesic all-set-homogeneous spaces include complete simply-connected Riemannian manifolds with constant curvature and the circle equipped with length metrics.
These will be referred further as \emph{classical spaces};
they are closely related to classical Euclidean/non-Euclidean geometry.

It is worth mentioning that an infinite-dimensional Hilbert space is \emph{not} all-set-homogeneous;
indeed, it is isometric to its proper subset.
Also, for $n\ge2$, the real projective space $\RP^n$ with canonical metric is not all-set-homogeneous.
Indeed, it contains two isometric but noncongruent triples of points with pairwise distance $\tfrac\pi3$
(we assume that a closed geodesic on $\RP^n$ has length $\pi$);
one triple lies on a closed geodesic and another does not. 

Nonclassical examples include the universal metric trees of finite valence;
these are discussed in the next section.

Given a metric space $M$ and a positive integer $n$, consider all pseudometrics induced on $n$ points $x_1,\dots, x_n\in M$.
Any such metric is completely described by $N=\tfrac{n\cdot (n-1)}2$ distances $|x_i-x_j|_M$ for $i<j$, so it can be encoded by a point in $\RR^N$.
The set of all these points $F_n(M)\subset \RR^N$ will be called \emph{$n^\text{th}$ fingerprint} of~$M$.

\begin{thm}{Theorem}\label{all-sets}
Let $M$ be a complete all-set-homogeneous length space.
Suppose that all fingerprints of $M$ are closed.
Then $M$ is classical.
\end{thm}

The following two results are closely related to our theorem.
\begin{itemize}
\item \emph{Any complete all-set-homogeneous geodesic space with locally unique nonbifurcating geodesics is classical;} it was proved by Garrett Birkhoff \cite{birkhoff}.
\item \emph{Any locally compact three-point-homogeneous length space is classical.}
 This result was proved by Herbert Busemann \cite{busemann}; it also follows from the more general result of Jacques Tits \cite{tits} about two-point-homogeneous spaces.
 (A space is called $n$-point homogeneous if any isometry between its subsets with at most $n$ points in each can be extended to an isometry of the whole space.)
\end{itemize}
For more related results, see the survey by Semeon Bogatyi \cite{bogaty} and the references therein.

\parit{Proof.}
If $M$ is locally compact, then the statement follows from the Busemann--Tits result stated above.
Therefore, we can assume that $M$ is not locally compact.

In this case, there is an infinite sequence of points $x_1,x_2,\dots$ such that 
$\eps<|x_i-x_j|_M<1$ for some fixed $\eps>0$ and all $i\ne j$.
Applying the Ramsey theorem, we get that for arbitrary positive integer $n$ there is a sequence $x_1,x_2,\dots,x_n$ such that all the distances
$|x_i-x_j|_M$ lie in an arbitrarily small subinterval of $(\eps,1)$.
Since the fingerprints are closed, there is an arbitrarily long sequence 
$x_1,x_2,\dots,x_n$ such that 
$|x_i\z-x_j|_M= r$ for some fixed $r>0$.

Choose a maximal (with respect to inclusion) set of points $A$ with distance $r$ between any pair.
Since $M$ is all-set-homogeneous, we get that $A$ has to be infinite.
In particular, there is a map $f\:A\to A$ that is injective, but not surjective.

Note that $f$ is distance-preserving.
Since $A$ is maximal, for any $y\notin A$ we have that $|y-x|_M\ne r$ for some $x\in A$.
It follows that a distance-preserving map $M\to M$ that agrees with $f$ cannot have points of $A\setminus f(A)$ in its image.
In particular, no isometry $M\to M$ agrees with the map $f$ --- a contradiction.
\qeds

\section{Example}

Recall that geodesic space $T$ is called a \emph{metric tree} if any pair of points $x,y\in T$ are connected by a unique geodesic $[xy]_T$,
and the union of any two geodesics $[xy]_T$, and $[yz]_T$ contains $[xz]_T$.
The \emph{valence} of $x\in T$ is defined as the cardinality of connected components in $T\setminus \{x\}$.

It is known that for any cardinality $n\ge 2$, there is a space $\TT_n$ that satisfies the following properties:
\begin{itemize}
\item The space $\TT_n$ is a complete metric tree with valence $n$ at any point.
\item $\TT_n$ is homogeneous; that is, its group of isometries acts transitively. 
\end{itemize}
Moreover, this space is uniquely defined up to isometry and $n$-universal; the latter means that $\TT_n$ includes an isometric copy of any metric tree of maximal valence at most $n$.

The space $\TT_n$ is called
a \emph{universal metric tree of valence $n$}.
An explicit construction of $\TT_n$ is given by Anna Dyubina and Iosif Polterovich~\cite{dyubina-polterovich}.
Their proof of the universality of $\TT_n$ admits a straightforward modification that proves the following claim. 

\begin{thm}{Claim}
If $n$ is finite, then $\TT_n$ is all-set-homogeneous.
\end{thm}

Note that the claim implies that the condition on fingerprints in the theorem is necessary.
In fact, if $n\ge 3$, then the $(n+1)^{\text{th}}$ fingerprint of $\TT_n$ is not closed --- $\TT_n$ does not contain $n+1$ points on distance 1 from each other,
but it contains an arbitrarily large set with pairwise distances arbitrarily close to 1.

\parit{Proof.}
Let $A, A'\subset \TT_n$ and $x\mapsto x'$ be an isometry $A\to A'$.
Applying the Zorn lemma, we can assume that $A$ is maximal; that is, the domain $A$ cannot be extended by a single point.
It remains to show that $A=\TT_n$ and $A'=\TT_n$.

Note that $A$ is closed.

\begin{wrapfigure}{o}{23 mm}
\vskip-6mm
\centering
\includegraphics{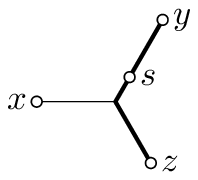}
\end{wrapfigure}

Further, suppose $x,y,z\in A$ and $s\in [yz]_{\TT_n}$.
Since $\TT_n$ is a metric tree, the distance $|x-s|_{\TT_n}$ is completely determined by four values $|x-y|_{\TT_n}$, $|x-z|_{\TT_n}$, $|s-y|_{\TT_n}$, $|s-z|_{\TT_n}$.

Denote by $s'$ the point on the geodesic $[y'z']_{\TT_n}$ such that $|y'-s'|_{\TT_n}\z=|y-s|_{\TT_n}$ and therefore $|z'-s'|_{\TT_n}=|z-s|_{\TT_n}$.
Since the map preserves distances $|x-y|_{\TT_n}$ and $|x-z|_{\TT_n}$, we get $|s'-x'|_{\TT_n}=|s-x|_{\TT_n}$;
that is, the extension of the map by $s\mapsto s'$ is still distance-preserving.

Since $A$ is maximal, $s\in A$.
In other words, $A$ is a convex subset of $\TT_n$;
in particular, $A$ is a metric tree with maximal valence at most~$n$.

Arguing by contradiction, suppose $A\ne \TT_n$, choose $a\in A$ and $b\notin A$. 
Let $c\in A$ be the last point on the geodesic $[ab]_{\TT_n}$.
Note that the valence of $c$ in $A$ is smaller than $n$.

Since $n$ is finite, at least one of the connected components in $\TT_n\z\setminus \{c'\}$ does not intersect $A'$.
Choose a point $b'$ in this component such that $|c'-b'|_{\TT_n}=|c-b|_{\TT_n}$.
Observe that the map can be extended by $b\mapsto b'$ --- a contradiction.
It follows that $A=\TT_n$.

It remains to show that $A'=\TT_n$.
Note that $A'$ is a closed convex set in $\TT_n$ that is isometric to $\TT_n$.
In particular valence of any point in $A'$ is $n$.

Assume $A'$ is a proper subset of $\TT_n$.
Choose $a'\in A'$ and $b'\notin A'$.
Let $c'\in A'$ be the last point on the geodesic $[a'b']_{\TT_n}$.
Observe that the valence of $c'$ in $A'$ is smaller than $n$ --- a contradiction.
\qeds

\section{Remarks}

Let us list examples for related classification problems.
We would like to see any other example or a proof of the corresponding classification. 

First of all, we do not see other examples of complete all-set-homogeneous length spaces except those listed in the theorem and the claim.

Without length-metric assumption, we have a vast amount of examples.
It includes finite discrete spaces, Cantor sets with natural ultrametrics;
also note that snowflaking $(X,|\ -\ |^\theta)$ of any all-set-homogeneous spaces $(X,|\ -\ |)$ is all-set-homogeneous.

The definition of all-set-homogeneous spaces can be restricted to \emph{small} subsets $A$ and $A'$; for example, \emph{finite} or \emph{compact}.
In these cases, we say that the space is \emph{finite-set-homogeneous} or \emph{compact-set-homogeneous} respectively.

Examples of complete separable compact-set-homogeneous length spaces include the spaces listed in the theorem,
plus the Urysohn spaces $\UU$ and $\UU_d$ (the space $\UU_d$ is isometric to a sphere of radius $\tfrac d2$ in $\UU$).
Without the separability condition, we get in addition the metric trees from the claim.

The finite-set-homogeneous spaces include, in addition, infinite-dimensional analogs of the classical spaces;
in particular the Hilbert space.  

Let us also mention that finite-set homogeneity is closely related to the \emph{metric version of Fraïssé limit} introduced by Itay Ben Yaacov \cite{ben-yaacov}. 

\parbf{Acknowledgments.}
This note is inspired by the question of Joseph O'Rourke \cite{rourke}.
We want to thank James Hanson for his interesting and detailed comments on our question \cite{hanson}.
The second author wants to thank Rostislav Matveyev for an interesting discussion on Rubinstein Street. 

The first author was partially supported by the Russian Foundation for Basic Research grant 20-01-00070; the second author was partially supported by the National Science Foundation grant DMS-2005279 and the Ministry of Education and Science of the Russian Federation, grant 075-15-2022-289.

{\sloppy
\printbibliography[heading=bibintoc]

@article {tits,
    AUTHOR = {Tits, J.},
     TITLE = {Sur certaines classes d'espaces homog\`enes de groupes de {L}ie},
   JOURNAL = {Acad. Roy. Belg. Cl. Sci. M\'{e}m. Coll. in 8$^\circ$},
  FJOURNAL = {Acad\'{e}mie Royale de Belgique. Classe des Sciences. M\'{e}moires.
              Collection in-8$^\circ$. Koninklijke Belgische Academie.
              Klasse der Wetenschappen. Verhandelingen. Verzameling
              in-8$^\circ$},
    VOLUME = {29},
      YEAR = {1955},
    NUMBER = {3},
     PAGES = {268},
      ISSN = {0365-0936},
   MRCLASS = {17.0X},
  MRNUMBER = {76286},
MRREVIEWER = {L. Auslander},
}

@article {ben-yaacov,
    AUTHOR = {Ben Yaacov, I.},
     TITLE = {Fraïssé limits of metric structures},
   JOURNAL = {J. Symb. Log.},
  FJOURNAL = {The Journal of Symbolic Logic},
    VOLUME = {80},
      YEAR = {2015},
    NUMBER = {1},
     PAGES = {100--115},
      ISSN = {0022-4812},
   MRCLASS = {03C30 (03C52 03C68)},
  MRNUMBER = {3320585},
MRREVIEWER = {M. Yasuhara},
       DOI = {10.1017/jsl.2014.71},
       URL = {https://doi.org/10.1017/jsl.2014.71},
}

@article {birkhoff,
    AUTHOR = {Birkhoff, G.},
     TITLE = {Metric foundations of geometry. {I}},
   JOURNAL = {Trans. Amer. Math. Soc.},
  FJOURNAL = {Transactions of the American Mathematical Society},
    VOLUME = {55},
      YEAR = {1944},
     PAGES = {465--492},
      ISSN = {0002-9947},
   MRCLASS = {48.0X},
  MRNUMBER = {10393},
MRREVIEWER = {L. M. Blumenthal},
       DOI = {10.2307/1990304},
       URL = {https://doi.org/10.2307/1990304},
}

@article {bogaty,
    AUTHOR = {Bogatyi, S. A.},
     TITLE = {Metrically homogeneous spaces},
   JOURNAL = {Russian Math. Surveys},
    VOLUME = {57},
      YEAR = {2002},
    NUMBER = {2},
     PAGES = {221–-240},
      ISSN = {0042-1316},
   MRCLASS = {54F65 (53C70)},
  MRNUMBER = {1918193},
MRREVIEWER = {Vitali Kapovitch},
       DOI = {10.1070/RM2002v057n02ABEH000495},
       URL = {https://doi.org/10.1070/RM2002v057n02ABEH000495},
}

@article {bogaty-ru,
    AUTHOR = {Богатый, С. А.},
     TITLE = {Метрически однородные пространства},
   JOURNAL = {Успехи математических наук},
    VOLUME = {57},
      YEAR = {2002},
    NUMBER = {2},
     PAGES = {3---22},
      ISSN = {0042-1316},
   MRCLASS = {54F65 (53C70)},
  MRNUMBER = {1918193},
MRREVIEWER = {Vitali Kapovitch},
       DOI = {10.1070/RM2002v057n02ABEH000495},
       URL = {https://doi.org/10.1070/RM2002v057n02ABEH000495},
       language={russian},
       hypernation={russian},
}

@book {busemann,
    AUTHOR = {Busemann, H.},
     TITLE = {Metric methods in {F}insler spaces and in the
              foundations of geometry},
    SERIES = {Annals of Mathematics Studies, No. 8},
 %PUBLISHER = {Princeton University Press, Princeton, N. J.},
      YEAR = {1942},
  %   PAGES = {viii+243},
   MRCLASS = {48.0X},
  MRNUMBER = {0007251},
MRREVIEWER = {S. M. Ulam},
}

@article {dyubina-polterovich,
    AUTHOR = {Dyubina, A. and Polterovich, I.},
     TITLE = {Explicit constructions of universal {$\mathbb{R}$}-trees and
              asymptotic geometry of hyperbolic spaces},
   JOURNAL = {Bull. London Math. Soc.},
  FJOURNAL = {The Bulletin of the London Mathematical Society},
    VOLUME = {33},
      YEAR = {2001},
    NUMBER = {6},
     PAGES = {727--734},
      ISSN = {0024-6093},
   MRCLASS = {57M07 (20F67 53C23 54F50)},
  MRNUMBER = {1853785},
       DOI = {10.1112/S002460930100844X},
       URL = {https://doi.org/10.1112/S002460930100844X},
}

@MISC {hanson,
    TITLE = {All-set-homogeneous spaces},
    AUTHOR = {J. Hanson},
    HOWPUBLISHED = {MathOverflow},
   % NOTE = {URL:https://mathoverflow.net/q/430738 (version: 2022-09-19)},
    EPRINT = {https://mathoverflow.net/q/430738},
   % URL = {https://mathoverflow.net/q/430738}
}

@MISC {rourke,
    TITLE = {Which metric spaces have this superposition property?},
    AUTHOR = {J. O'Rourke},
    HOWPUBLISHED = {MathOverflow},
   % NOTE = {URL:https://mathoverflow.net/q/118008 (version: 2014-09-25)},
    EPRINT = {https://mathoverflow.net/q/118008},
   % URL = {https://mathoverflow.net/q/118008}
}
\fussy
}

\end{document}


\title{\thetitle}
\author{\theauthors}
\date{}
\maketitle

\begin{abstract}
Метрическое пространство называется полностью однородным, если любая изометрия между его подмножествами продолжается до изометрии всего пространства.
Мы классифицируем определённый подкласс полностью однородных пространств с внутренней метрикой.
\end{abstract}

\section{Основной результат}

Расстояние между точками $x$ и $y$ метрического пространства $M$ будет обозначаться как $|x-y|_M$.
Метрика на $M$ называется \emph{внутренней} (\emph{геодезической}) если любые две точки $x,y\in M$ можно соединить путём $\gamma$ таким, что расстояние $|x-y|_M$ произвольно близко к длине $\gamma$ (или, соответственно, равно этой длине).
Очевидно, что любая геодезическая метрика внутренняя, но не наоборот.

Метрическое пространство $M$ называется \emph{полностью однородным}, если любая изометрия $A\to A'$ между его подмножествами продолжается до изометрии $M\to M$ всего пространства.

Примеры геодезических полностью однородных пространств включают односвязные римановы многообразия постоянной кривизны и окружности с внутренней метрикой.
Эти примеры будут называться \emph{классическими пространствами},
они тесно связанны с классической евклидовой/неевклидовой геометрией.

Стоит заметить, что бесконечномерное гильбертово пространство не полностью однородно --- это следует из того, что оно изометрично своему собственному подмножеству.
Далее, при $n\ge2$, вещественное проективное пространство $\RP^n$ с канонической метрикой не полностью однородно.
Действительно, $\RP^n$ содержит две изометричные, но не конгруэнтные тройки точек с попарными расстояниями $\tfrac\pi3$ ---
одна тройка лежит на замкнутой геодезической, а другая нет. 
(Мы предполагаем, что длина замкнутой геодезической на $\RP^n$ равна $\pi$.)

Неклассические примеры включают универсальные деревья конечной валентности, они обсуждаются в следующем разделе.

Для метрического пространства $M$ и целого числа $n\ge1$ рассмотрим все псевдометрики на $n$-ках точек 
$x_1,\dots, x_n\in M$.
Каждая такая метрика описывается $N=\tfrac{n\cdot (n-1)}2$ расстояниями $|x_i-x_j|_M$ для $i<j$, и, таким образом, соответствует точке в $\RR^N$.
Множество всех таких точек $F_n(M)\z\subset \RR^N$ будет называться $n$-м отпечатком $M$.

\begin{thm}{Теорема}\label{all-sets}
Пусть $M$ --- полное полностью однородное пространство с внутренней метрикой.
Предположим, что все отпечатки $M$ замкнуты.
Тогда $M$ --- классическое пространство.
\end{thm}

Следующие два утверждения очень близки нашему результату:
\begin{itemize}
\item \emph{Любое полностью однородное геодезическое пространство с локально однозначными не расщепляющимися геодезическими является классическим.} Доказано Гарретом Биркгофом  \cite{birkhoff}.
\item \emph{Любое локально компактное 3-точечно однородное пространство с внутренней метрикой  является классическим.} 
Доказано Гербертом Буземаном \cite{busemann}.
Это утверждение также следует из более общего результата Жака Титса \cite{tits} о двуточечно однородных пространствах.
(Пространство называется $n$-точечно однородным, если любая изометрия между его подмножествами с не более чем $n$ точками в каждом продолжается до изометрии всего пространства.)
\end{itemize}
Подробное обсуждение этих и других связанных результатов дано в обзоре Семёна Богатого \cite{bogaty-ru}.

\parit{Доказательство.}
Если $M$ локально компактно, то утверждение следует из приведённой выше теоремы Буземана --- Титса.
Поэтому можно предположить, что $M$ не локально компактно.

В этом случае существует бесконечная последовательность точек $x_1$, $x_2,\dots$ такая, что $\eps<|x_i-x_j|_M<1$ для некоторого $\eps>0$ и всех $i\ne j$.
Применив теорему Рамсея, получаем, что для произвольного целого $n\ge 1$ существует последовательность $x_1,x_2,\dots,x_n$ такая, что все расстояния $|x_i-x_j|_M$ лежат в произвольно малом подинтервале $(\eps,1)$.
Поскольку отпечатки замкнуты, существует произвольно длинная последовательность $x_1,x_2,\dots,x_n$ такая, что
$|x_i\z-x_j|_M= r$ для некоторого фиксированного $r>0$ и всех $i\ne j$.

Выберем максимальное (относительно включения) множество $A$ с расстоянием $r$ между любыми двумя его точками.
Поскольку $M$ полностью однородно, множество $A$ обязано быть бесконечным.
В частности, существует отображение $f\:A\to A$, которое инъективно, но не сюрьективно.

Заметим, что $f$ сохраняет расстояния.
Поскольку $A$ максимально, для любой точки $y\notin A$ существует $x\in A$ такая, что $|y-x|_M\ne r$.
Значит любое отображение $M\to M$, которое согласуется с $f$ и сохраняет расстояния, не может содержать в своём образе точек из $A\setminus f(A)$.
В частности, никакая изометрия $M\to M$ не согласуется с $f$ --- противоречие.
\qeds

\section{Пример}

Напомним, что пространство $T$ с геодезической метрикой называется \emph{метрическим деревом} если любая пара точек $x,y\in T$ соединяется единственной кратчайшей $[xy]_T$ и объединение любой пары кратчайших $[xy]_T$ и $[yz]_T$ содержит $[xz]_T$.
\emph{Валентностью} точки $x\in T$ называется кардинальное число связных компонент в $T\setminus\{x\}$.

Для любого кардинального числа $n\ge 2$ сущестует пространство $\TT_n$ удовлетворяющее следующим условиям:
\begin{itemize}
\item $\TT_n$ --- полное метрическое дерево с валентностью $n$ в каждой точке.
\item $\TT_n$ --- однородно, то есть, группа его изометрий действует транзитивно. 
\end{itemize}
Более того, $\TT_n$ однозначно определено с точностью до изометрии и является $n$-универсальным.
Последнее означает, что $\TT_n$ включает в себя изометрическую копию любого метрического дерева с максимальной валентностью не больше $n$.

Пространство $\TT_n$ называется \emph{универсальным метрическим деревом валентности $n$}.
Явное построение $\TT_n$ дано Анной Дюбиной и Иосифом Полтеровичем~\cite{dyubina-polterovich}.
Доказательство следующего утверждения --- простая переделка их доказательства универсальности $\TT_n$.

\begin{thm}{Утверждение}
Если $n$ конечно, то $\TT_n$ полностью однородно.
\end{thm}

В частности, условие на отпечатки в нашей теореме является существенным.
Заметим, что если $n\ge 3$, то $(n+1)$-й отпечаток $\TT_n$ не замкнут --- $\TT_n$ не содержит $(n+1)$-точечное множество с попарными расстояниями 1, 
но содержит произвольно большое конечное множество с  попарными расстояниями произвольно близкими к 1.

\parit{Доказательство.}
Пусть $A, A'\subset \TT_n$, и $x\mapsto x'$ есть изометрия $A\to A'$.
По лемме Цорна, можно предположить, что $A$ максимально,
то есть, к $A$ невозможно добавить точку.
Остаётся доказать, что $A=\TT_n$ и $A'=\TT_n$.

Заметим, что $A$ замкнуто.

\begin{wrapfigure}{o}{23 mm}
\vskip-6mm
\centering
\includegraphics{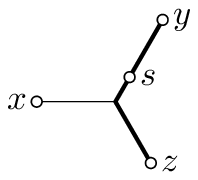}
\end{wrapfigure}

Далее, предположим $x,y,z\in A$ и $s\in [yz]_{\TT_n}$.
Поскольку $\TT_n$ --- метрическое дерево, расстояние $|x-s|_{\TT_n}$ полностью определяется четвёркой чисел 
$|x-y|_{\TT_n}$, $|x-z|_{\TT_n}$, $|s-y|_{\TT_n}$, $|s-z|_{\TT_n}$.

Обозначим через $s'$ точку на кратчайшей $[y'z']_{\TT_n}$ такую, что $|y'-s'|_{\TT_n}\z=|y-s|_{\TT_n}$, и значит $|z'-s'|_{\TT_n}=|z-s|_{\TT_n}$.
Поскольку отображение сохраняет расстояния $|x-y|_{\TT_n}$ и $|x-z|_{\TT_n}$,
получаем, что $|s'-x'|_{\TT_n}=|s-x|_{\TT_n}$.
То есть отображение продолженное как $s\mapsto s'$ остаётся сохраняющим расстояния.

Поскольку $A$ максимально, $s\in A$.
Другими словами, $A$ --- выпуклое подмножество в $\TT_n$.
В частности, $A$ --- метрическое дерево максимальной валентности не больше $n$.

Предположим $A\ne \TT_n$.
Выберем $a\in A$ и $b\notin A$.
Пусть $c\in A$ --- последняя точка на кратчайшей $[ab]_{\TT_n}$.
Заметим, что валентность $c$ в $A$ строго меньше $n$.

Поскольку $n$ конечно, хотя бы одна из компонент связности $\TT_n\z\setminus \{c'\}$ не пересекает $A'$.
Выберем точку $b'$ в этой компоненте так, чтобы $|c'-b'|_{\TT_n}\z=|c-b|_{\TT_n}$.
Заметим, что изометрию $A\to A'$ можно продолжить как $b\mapsto b'$ --- противоречие.
Значит $A=\TT_n$.

Осталось показать, что $A'=\TT_n$.
Заметим, что $A'$ замкнутое подмножество $\TT_n$ и оно изометрично $\TT_n$.
В частности валентность любой точки $A'$ равна $n$.

Предположим $A'$ --- собственное подмножество $\TT_n$.
Выберем $a'\in A'$ и $b'\notin A'$.
Пусть $c'\in A'$ --- последняя точка на кратчайшей $[a'b']_{\TT_n}$.
Тогда валентность $c'$ в $A'$ меньше $n$ --- противоречие.
\qeds

\section{Замечания}

Приведём известные нам примеры пространств для похожих классификационных задач.
Нас интересуют любые другие примеры а также доказательства таких классификаций.

Во-первых, мы не знаем примеров полных полностью однородных пространств с внутренней метрикой кроме классических пространств и универсальных метрических деревьев конечной валентности.

Без условия внутренности метрики, дополнительных примеров становится слишком много.
Они включают в себя конечные дискретные пространства, канторовы множества с естественными ультраметриками.
Также заметим, что \emph{снежинка} $(X,|\ -\ |^\theta)$ любого полностью однородного пространства $(X,|\ -\ |)$ полностью однородна.

В определении полностью однородных пространств можно рассматривать только маленькие подмножества $A$ и $A'$,
например, \emph{конечные} или \emph{компактные}.
Тогда мы говорим, что пространство  \emph{конечно однородно} или  \emph{компактно однородно} соответственно.

В дополнение к классическим пространствам,
примеры полных сепарабельных компактно однородных пространств с внутренней метрикой включают пространства Урысона  $\UU$ и $\UU_d$ (пространство $\UU_d$ изометрично сфере радиуса $\tfrac d2$ в $\UU$).
Без сепарабельности к ним добавляются универсальные метрические деревья конечной валентности.

Конечно однородные пространства включают дополнительно бесконечномерные аналоги классических пространств,
в частности гильбертово пространство.

Заметим также, что конечно однородные пространства тесно связаны с \emph{метрическим пределом Фраисэ} предложенным Итайем Бен Яаковом \cite{ben-yaacov}. 

\parbf{Благодарности.}
Эта заметка появилась благодаря вопросу Джозефа Орурка \cite{rourke}.
Мы хотим поблагодарить Джеймса Хансона за его интересные и подробные замечания к нашему вопросу \cite{hanson}.
Второй автор признателен Ростиславу Матвееву за интересный разговор на улице Рубинштейна.

Первый автор частично поддержан грантом РФФИ 20-01-00070;
второй автор частично поддержан грантом NSF DMS-2005279
и грантом Минобрнауки 075-15-2022-289.

{\sloppy
\printbibliography[heading=bibintoc]
\fussy
}